\documentclass[11pt]{article}

\usepackage{amsfonts,amsmath,latexsym,color,epsfig,hyperref,enumitem, amssymb}
\newtheorem{theorem}{Theorem}
\newtheorem{lemma}{Lemma}
\newtheorem{proposition}{Proposition}

\newtheorem{claim}{Claim}

\newcommand{\definition}
      {\medskip\noindent {\bf Definition:\hspace{0em}}}

\newenvironment{proof}
      {\medskip\noindent{\bf Proof:}\hspace{1mm}}
      {\hfill$\Box$\medskip}
\def\qed{\ifvmode\mbox{ }\else\unskip\fi\hskip 1em plus 10fill$\Box$}
\newenvironment{proofof}[1]
      {\medskip\noindent{\bf Proof of #1:}\hspace{1mm}}
      {\hfill$\Box$\medskip}
\input{epsf}

\makeatletter
\def\Ddots{\mathinner{\mkern1mu\raise\p@
\vbox{\kern7\p@\hbox{.}}\mkern2mu
\raise4\p@\hbox{.}\mkern2mu\raise7\p@\hbox{.}\mkern1mu}}
\makeatother

\newcommand{\F}{\mathbb F}
\newcommand{\N}{\mathbb N}

\title{\vspace{-0.7cm}Zero subsums in vector spaces over finite fields}
\author{Cosmin Pohoata\thanks{Department of Mathematics, Yale University, USA. Email: {\tt andrei.pohoata@yale.edu}.} \and Dmitriy Zakharov\thanks{Laboratory of Combinatorial and Geometric Structures, MIPT, Russia. Email: {\tt  zakharov2k@gmail.com}.}}
\date{}

\begin{document}
\maketitle

\begin{abstract}
The Olson constant $\mathcal{O}L(\mathbb{F}_{p}^{d})$ represents the minimum positive integer $t$ with the property that every subset $A\subset \mathbb{F}_{p}^{d}$ of cardinality $t$ contains a nonempty subset with vanishing sum. The problem of estimating $\mathcal{O}L(\mathbb{F}_{p}^{d})$ is one of the oldest questions in additive combinatorics, with a long and interesting history even for the case $d=1$.

In this paper, we prove that for any fixed $d \geq 2$ and $\epsilon > 0$, the Olson constant of $\mathbb{F}_{p}^{d}$ satisfies the inequality
$$\mathcal{O}L(\mathbb{F}_{p}^{d}) \leq (d-1+\epsilon)p$$
for all sufficiently large primes $p$. This settles a conjecture of Hoi Nguyen and Van Vu.
\end{abstract}

\section{Introduction}

\smallskip

For a subset $A$ of an additive group $G$, consider the set of all nonempty subsums
$$\Sigma^*(A):=\left\{\sum_{x\in B}x\ |\ B\subset A, B\ne \emptyset\right\}.$$
The \emph{Olson constant} $\mathcal{O}L(G)$ represents the minimum $t$ such that every subset $A\subset G$ of cardinality $t$ satisfies $0\in \Sigma^*(A)$. This is a well-known quantity in additive combinatorics, which is notoriously difficult to estimate even for the most basic groups. Its nice history begins in 1964 with Erd\H{o}s and Heilbronn, who in \cite{EH64} proved that there exists an absolute constant $c$ such that $\mathcal{O}L(\mathbb{F}_{p}) \leq c\sqrt{p}$, where $p$ is an odd prime. In the same paper, they conjectured that their result should generalize to arbitrary additive groups $G$ and that the optimal constant $c$ in the inequality above is probably $c = \sqrt{2}$. A few years later Szemer\'edi \cite{Sze70} settled the former conjecture in the affirmative. The result of Erd\H{o}s and Heilbronn for $\mathbb{F}_{p}$ and Szemer\'edi's theorem for general groups were both later refined by Olson in \cite{Ols69}, \cite{Ols692} and \cite{Ols75}, who proved that $\mathcal{O}L(G) \leq 2\sqrt{|G|}$ and also introduced a remarkable group ring approach (which has also recently resurfaced in the context of the polynomial method developments around the cap set problem; see \cite{P17} and \cite{PP20}). Olson's result was subsequently pushed further by Hamidoune and Zemor \cite{HZ96}, who proved that $\mathcal{O}L(G) \leq \sqrt{2|G|} + O(|G|^{1/3} \log|G|)$, and among other things established the correct order of growth for $\mathcal{O}L(\mathbb{F}_{p})$, up to lower order terms. In 2008, Nguyen, Szemer\'edi and Vu \cite{NSV08} finally removed the lower terms in the primordial case $G = \mathbb{F}_{p}$, therefore proving the optimal inequality $\mathcal{O}L(\mathbb{F}_{p}) \leq \sqrt{2p}$ for all sufficiently large primes $p$. This work was also further refined in two separate rounds by Balandraud in \cite{Bal12} and \cite{Bal17}, who finally gave a short alternative argument which works for {\textit{all}} odd primes $p$ based on the quantitative Combinatorial Nullstellensatz introduced by Karasev and Petrov in \cite{KP12}.

In this paper, we address the problem for $G=\mathbb{F}_{p}^{d}$, where $p$ is an odd prime number, $d \geq 2$, and $\mathbb{F}_{p}^{d}$ denotes as usual the vector space of $d$-dimensional vectors with coordinates from $\mathbb{F}_{p}$. The situation in higher dimensions has been traditionally known to be much more complicated. For $d=2$, the first important result only appeared in 2004. In \cite{GRT04}, Gao, Ruzsa and Thangadurai proved that $\mathcal{O}L(\mathbb{F}_{p}^{2})=p+\mathcal{O}L(\mathbb{F}_{p})-1$ holds for all primes $p > 4.67 \times 10^{34}$, thus establishing a beautiful connection between $\mathcal{O}L(\mathbb{F}_{p}^{2})$ and $\mathcal{O}L(\mathbb{F}_{p})$. In particular, given the successful story for $\mathbb{F}_{p}$, this result also determines the Olson constant constant of $\mathbb{F}_{p}^{2}$ for large primes. For higher dimensions, however, not much more is known. In the same paper \cite{GRT04}, Gao, Ruzsa and Thangadurai conjectured that 
\begin{equation} \label{GRT}
\mathcal{O}L(\mathbb{F}_{p}^{d}) = p + \mathcal{O}L(\mathbb{F}_{p}^{d-1})-1
\end{equation}
should hold in general for all $d \geq 2$ and for all sufficiently large primes $p$, but this is still a (difficult) open problem. It is also perhaps worth mentioning the curiosity that the assumption that $p$ is sufficiently large is necessary this time around, see for instance the discussion from \cite{GG06}. In 2011, Nguyen and Vu \cite{NV12} also studied this higher dimensional problem and proposed the following asymptotic version of the conjecture: for any fixed $d \geq 2$ and $\epsilon > 0$, the Olson constant of $\mathbb{F}_{p}^{d}$ satisfies the inequality
\begin{equation} \label{Vu}
\mathcal{O}L(\mathbb{F}_{p}^{d}) \leq (d-1+\epsilon)p
\end{equation}
for all sufficiently large primes $p$. Since $\mathcal{O}L(\mathbb{F}_{p}) = O(\sqrt{p}) = o(p)$, it is clear that \eqref{GRT} implies \eqref{Vu}, but in some sense \eqref{Vu} should still capture all of the difficulties around \eqref{GRT} when $d \geq 3$. Extending an elegant alternative approach they found for the case $d=2$ of the Gao-Ruzsa-Thangadurai conjecture, Nguyen and Vu then established this asymptotic conjecture when $d=3$; however, their argument has various serious limitations already starting with $d \geq 4$, and so no further progress has been made since.

Our main result is a resolution of this conjecture of Nguyen and Vu in all dimensions $d \geq 2$ by using a new approach inspired by the second author's recent work on the Erd\H{o}s-Ginzburg-Ziv problem \cite{Zak20}. 

\begin{theorem} \label{main}
For any fixed $d \geq 2$ and $\epsilon > 0$, the Olson constant of $\mathbb{F}_{p}^{d}$ satisfies the inequality
$$\mathcal{O}L(\mathbb{F}_{p}^{d}) \leq (d-1+   \epsilon)p$$
for all sufficiently large primes $p$. 
\end{theorem}

We include the proof of Theorem \ref{main} in Section 3, after discussing terminology and the required preliminary results in Section 2. 

Before we move on however to the technical details, we end this section with a high-level overview of the argument. Starting with a set $X \subset \F_p^d$ of size $(d-1+\varepsilon)p$, where $p$ is a sufficiently large prime number, the first important idea is to prove that one can reduce the problem of finding a vanishing subsum in $\Sigma^{*}(X)$ to the case when $X$ lies in a translate of the form $v+[-K, K]^l \times \F_p^{d-l}$, for some $l \in \left\{1,\ldots,d-1\right\}$, $v \in \F_p^d$, and where $[-K,K]$ stands for the interval $\left\{-K,-(K-1),\ldots,(K-1),K\right\}$ --regarded as a subset of $\mathbb{F}_{p}$ (whose size does not depend on $p$). The second idea is that if $0 \not\in \Sigma^{*}(X)$, then one can also force $X$ to satisfy some further refined structural properties. Roughly, we'll be able to assume among other things, for example, that $X$ must always have some positive proportion of its elements outside the set $\left\{ x \in \F_p^d ~|~ \xi(x) \in [-K, K]\right\}$, for every linear function $\xi\ :\ \mathbb{F}_{p}^{d} \to \mathbb{F}_{p}$ (except for some trivial cases). The absence of ``linear concentration" is crucial because the third main idea is to consider the projection of this structured set $X$ onto the first $l$ coordinates. The image $Y$ of this projection is a large multiset in $\mathbb{F}_{p}^{l}$, so we can make use of tools such as the Combinatorial Nullstellensatz to find a suitable subsequence whose sum of elements vanishes and whose elements have various prescribed multiplicities. Finally, in order to close the argument, we then need to use the rich structure of $X$ to lift this auxiliary zero sum subsequence in $Y$ from the previous step up to an actual proper subset of $X$ whose sum of elements vanishes. 

While sharing a rather similar philosophy with the method of Nguyen and Vu from \cite{NV12} (where projection is also important), finding the right framework to project and lift in order to capture higher dimensional (additive) information and establishing the precise structural results which allow our procedure to go through for all dimensions $d \geq 2$ requires several new ideas, with both algebraic and probabilistic ingredients.

\bigskip

\section{Preliminaries}

\smallskip

A function $\xi$ on the space $\F_p^d$ is called {\textit{linear}} if it has the form 
\begin{equation}
    \xi(x_1, \ldots, x_d) = a_0 + a_1 x_1 +\ldots + a_d x_d,\label{lin}
\end{equation}
for some $a_i \in \F_p$. Linear functions $\xi_1, \ldots, \xi_l$ are called linearly independent if their ``linear parts", i.e. the vectors $(a_1, \ldots, a_d)$ from (\ref{lin}), are linearly independent in $\mathbb{F}_{p}^{d}$.

If $\xi$ is a linear function and $K \in \N$ we let
$$
H(\xi, K) :=\left\{ x \in \F_p^d ~|~ \xi(x) \in [-K, K]\right\},
$$
where $[-K,K]$ stands for the interval $\left\{-K,-(K-1),\ldots,(K-1),K\right\}$, regarded as a subset of $\mathbb{F}_{p}$. Given a linear function $\xi\ :\ \mathbb{F}_{p}^{d} \to \mathbb{F}_p$ and a multiset $X$ in $\F_p^d$, for $K \in \N$ and $\delta \ge 0$ we say that $X$ is {\textit{$(K, \delta)$-thick along $\xi$}} if $|X \setminus H(\xi, K)| \ge \delta |X|$, where the cardinality $|X|$ is calculated with multiplicities.

Sometimes it will be convenient to use functions $f: \mathbb F_p^d \rightarrow \mathbb R_+$ instead of multisets $X \subset \mathbb F_p^d$. For a subset $S \subset \mathbb F_p^d$ we denote by $f(S)$ the sum of values of $f$ on elements of $S$. We say that a function $f$ is {\it $(K, \delta)$-thick along a linear function $\xi$} if $f(H(\xi, K)) \le (1 - \delta) f(\mathbb F_p^d)$.

If $A$ is a multiset then {\it a submultiset (or subset for shortcut)} of $A$ is a multiset $B$ such that the multiplicity of every element in $B$ is at most the corresponding multiplicity in $A$. For technical reasons, we will also need a slightly different notion. A {labeled subset} of a multiset $A$ is a set $B \subset A \times \N$ such that every element $(a, n) \in B$ has the property that the multiplicity of $a$ in $A$ is at least $n$. So every finite multiset $A$ has exactly $2^{|A|}$ labeled subsets. Abusing notation, we sometimes also denote by $B$ the corresponding submultiset of $A$. A pair $(a, n) \in A \times \N$ where $n$ is at most the multiplicity of $a$ in $A$ is called {\it a labeled element of $A$}.

In what follows, $\delta \in (0, 1)$, $K'$ and $K$ are natural numbers, and for $n \in \N$ we denote by $[n]$ the set $\{1, \ldots, n\}$.

\smallskip

\subsection{Tube decomposition}

\smallskip

\definition{
A set $X \subset \F_p^d$ is called {\textit{$(K, K', \delta)$-tubular}} if there exists $l \in [0, d]$ and an affine isomorphism $\psi: \F_p^d \rightarrow \F_p^d$ such that $\psi X$ is contained in the set $[-K, K]^l \times \F_p^{d-l}$ and satisfies the following property: $\psi X$ is $(K', \delta)$-thick along any linear function $\xi$ which is not constant on $\{0\} \times \F_p^{d-l}$. 
}


\smallskip

So, for instance, if $l = 0$ then $X$ is $(K', \delta)$-thick along any non-constant linear function. If $l = d$ then, after an appropriate change of coordinates, $X$ is contained in the box $[-K, K]^d \subset \F_p^d$. In general, we allow some combination of the above situations.

Note that a set $X$ is $(0, K, \delta)$-tubular if it is $(K, \delta)$-thick along any linear map $\xi$ which is not constant on the affine hull of $X$. Indeed, the space $\{0\} \times \F_p^{d-l}$ from the definition above must coincide with the affine hull of $\psi X$.
We shall say that in this case $X$ is {\textit{$(K,\delta)$-thick in its affine hull}}.

\begin{lemma}\label{tube}
Let $g: \N \rightarrow \N$ be any increasing function and let $K_{0} \in \mathbb{N}$, $\delta \in (0,1)$, and $d \in \mathbb{N}$ be such that $\delta < 2^{-d-1}$. For any set $X \subset \F_p^d$, there exists $Y \subset X$ of size at least $(1 - 2^{d+1}\delta)|X|$ which is $(K, g(K), \delta)$-tubular where $K = g^l(K_0)$ for some $l \in [0, d]$. 
\end{lemma}


\begin{proofof}{Lemma \ref{tube}}
Let $\xi_1, \ldots, \xi_l$ be a maximal sequence of linearly independent linear functions such that $X$ is not $(g^i(K_0), 2^i \delta)$-thick along $\xi_i$ for any $i  = 1, \ldots, l$. Note that we can have $l = 0$, in which case we put $g^0(K_0) = K_0$.


Consider $K = g^l(K_0)$ and let 
$$
Y = X \cap \bigcap_{i = 1}^l H(\xi_i, K).
$$
By the definition of the $\xi_i$'s we have 
\begin{equation}\label{lb}
|Y| \ge |X| - |X|\sum_{i=1}^l 2^i \delta \ge |X| (1 - (2^{l+1}-1)\delta).
\end{equation}
Moreover, for any linear function $\eta$ which is linearly independent from $\xi_1, \ldots, \xi_l$, note that the set $X$ is $(g(K), 2^{l+1}\delta)$-thick along $\eta$. Consequently, by (\ref{lb}), the set $Y \subset X$ is $(g(K), \delta)$-thick along $\eta$. Moreover, after an appropriate change of coordinates, we have $Y \subset [-K, K]^l \times \F_p^{d-l}$, where the copy of $\mathbb{F}_{p}^{d-l}$ arises as the ($d-l$)-dimensional intersection of the kernels of the maps $\xi_1,\ldots,\xi_l$.

\end{proofof}

\begin{lemma}\label{dec}
For any increasing function $g: \N\rightarrow \N$, $K_0 \ge 0$, $\varepsilon > 0$, and $d\ge 1$ 
there is some $N = N(K_0, d, \varepsilon, g) \in \N$ and functions $\mu(\cdot, \varepsilon, d),\delta(\cdot, \varepsilon, d) :\ \mathbb{N}\to \mathbb{R}_{>0}$ such that the following holds:
For any (multi-)set $X \subset \F_p^d$ there is $l \in [0, N]$ and a decomposition 
\begin{equation}\label{deq}
X = X_0 \cup X_1 \cup \ldots \cup X_m,
\end{equation}
such that $|X_0| \le \varepsilon |X|$ and  for any $i \in [m]$ we have $|X_i| \ge \mu |X|$ and $X_i$ is $(g(K), \delta)$-thick in its affine hull. Here $K = g^{l}(K_0)$, $\mu = \mu(K, \varepsilon, d)$ and $\delta = \delta(K, \varepsilon, d)$.
\end{lemma}

\begin{proofof}{Lemma \ref{dec}}
The proof is by induction on $d$. Take an arbitrary (multi-)set $X \subset \F_p^d$. If $X$ is $(g(K_0), \varepsilon/2)$-thick along any non-constant linear map, then there is nothing to prove since we can take the decomposition $X_0 = \emptyset$, $X_1 = X$ and $l = 0$. So we may assume that $X$ is not $(g(K_0), \varepsilon/2)$-thick. Then after a change of coordinates and removing at most $\varepsilon/2 |X|$ elements from $X$ we may assume that $X \subset [-g(K_0), g(K_0)] \times \F_p^{d-1}$. For each $y \in [-g(K_0), g(K_0)]$ let $X_y = X \cap \left(\{y\} \times \F_p^{d-1}\right)$. Remove from $X$ all sets $X_y$ such that $|X_y| < \varepsilon|X|/8g(K_0)$, so that the size of $X$ will decrease at most by $\varepsilon |X|/4$. Denote $\varepsilon' = \varepsilon / 8g(K_0)$.

Now we are going to apply the induction hypothesis to each of the remaining sets $X_y$. Let $X_{y_1}, \ldots, X_{y_r}$ be the list of all these sets, where $r \le 2g(K_0)$. We apply induction to $X_{y_1}$ with $\varepsilon = \varepsilon'$, $K'_0 = g(K_0)$ and $g = g^{N_1}$ where $N_1$ will be determined later. In order to the number $K = g^{l_1 N_1}(K'_0)$ from the induction hypothesis to lie in the interval $[K_0, g^N(K_0)]$ we need the following inequality:
\begin{equation}\label{ch}
    N >  N(K'_0, d-1, \varepsilon', g^{N_1}) \cdot N_1.
\end{equation}
So we get a decomposition of the form $X_{y_1} = \bigcup_{i = 0}^{m_1} X_{1, i}$ and there is some $l_1 \le N(K'_0, d-1, \varepsilon', g^{N_1})$ so that if we let $K_1 = g^{l_1 N_1}(K'_0)$ then
$X_{1, i}$ is $(g^{N_1}(K_1), \delta_1)$-thick in its affine hull for every $i \in [m_1]$. 

Now we apply the induction hypothesis to the set $X_{y_2}$ with $\varepsilon = \varepsilon'$, $K_0 = K_1$ and $g = g^{N_2}$ where $N_2$ will be determined later. To apply induction we need the following inequality:
$$
N_1 > N(K_1, d-1, \varepsilon', g^{N_2}) \cdot N_2.
$$
We thus will obtain a decomposition $X_{y_2} = \bigcup_{i = 0}^{m_2} X_{2, i}$ where $X_{2, i}$ is $(g^{N_2}(K_2), \delta_2)$-thick in its affine hull where $K_2 = g^{N_2 l_2}(K_1)$ for some $l_2 \le N(K_1, d-1, \varepsilon', g^{N_2})$. Moreover, we have $\delta_2 \gg_{K_2, \varepsilon/8K_0, d} 1$. 

Now observe that we have the following chain of inequalities:
$$
K_1 \le K_2 \le g^{N_2}(K_2) \le g^{N_1}(K_1).
$$
Thus, for every $i \in [m_1]$, the sets $X_{1, i}$ is $(g^{N_2}(K_2), \delta_1)$-thick in its affine hull. Also note that since $K_2 \ge K_1$ we have $\delta_1 \gg_{K_2, \varepsilon', d} 1$ as well. 

Applying induction in a similar manner $r-2$ more times to sets $X_{y_j}$ for $j = 3, \ldots, r$ we will eventually get some $K_r = g^l(K_0)$ where $l \le N$ such that all sets $X_{j, i}$ are $(g(K_r), \delta_r)$-thick in their affine hulls for some $\delta_r \gg_{K_r, \varepsilon', d} 1$. We will get a chain of inequalities of the form (\ref{ch}) which will give an upper bound on the function $N(K_0, d, \varepsilon, g)$.
This concludes the proof.

\end{proofof}

We will need a stronger version of Lemma \ref{dec}:

\begin{lemma}\label{sdec}
For any increasing function $g: \N\rightarrow \N$, $K_0 \ge 0$, $\varepsilon > 0$, and $d\ge 1$ 
there is some $N = N(K_0, d, \varepsilon, g) \in \N$ and functions $\mu(\cdot, \varepsilon, d),\delta(\cdot, \varepsilon, d) :\ \mathbb{N}\to \mathbb{R}_{>0}$ such that the following holds:
For any set $X \subset \F_p^d$ there is $l \in [0, N]$ and a decomposition 
\begin{equation*}
X = X_0 \cup X_1 \cup \ldots \cup X_m,
\end{equation*}
such that $|X_0| \le \varepsilon |X|$ and  for any $i \in [m]$ we have $|X_i| \ge \mu |X|$ and $X_i$ is $(g^{d+1}(K), \delta)$-thick in its affine hull. Here $K = g^{l}(K_0)$, $\mu = \mu(K, \varepsilon, d)$ and $\delta = \delta(K, \varepsilon, d)$.

Moreover, for any $S \subset [m]$ the set $X_S = \bigcup_{i \in S} X_i$ is $(K_S, g(K_S), \mu)$-tubular where $K_S = g^s(K)$ for some $s \in [0, d]$.

\end{lemma}

\begin{proofof}{Lemma \ref{sdec}}
Let $g' = g^{d+1}$ and apply Lemma \ref{dec} to $X$ with $g'$ instead of $g$, $\varepsilon/2$ instead of $\varepsilon$ and $K_0=K_0$. We get a decomposition of the form (\ref{deq}) where sets $X_i$ are $(g^{d+1}(K), \delta_0)$-thick in their affine hulls for some $\delta_0 \gg_{K, \varepsilon, d} 1$. Also we have $|X_i| \ge \mu_0 |X|$ for some $\mu_0 \gg_{K, \varepsilon, d} 1$. 

Let $S_1, S_2, \ldots, S_{2^m-1}$ be the list of all non-empty subsets of $[m]$ in any order. For $j=1, \ldots, 2^m-1$, apply Lemma \ref{tube} consecutively to sets $\bigcup_{i \in S} X_i$ with $K_0 = K$, $g=g$ and 
$$
\delta_j = \varepsilon \mu_0 \delta_0 2^{-d-2 - m} 2^{-(d+m+4) j}.  
$$
Such choice of $\delta_j$ will guarantee us the following properties:
\begin{enumerate}
    \item The total number $R$ of removed elements will be at most 
$$
R \le 2^{d+1}|X|\sum_{j = 1}^{2^m-1} \delta_j < \varepsilon|X|/2.
$$
    \item All sets $X_i$ will be $(g^{d+1}(K), \delta_0/2)$-thick since for every $i \in [m]$:
$$
R \le 2^{d+1}|X|\sum_{j = 1}^{2^m-1} \delta_j \le \varepsilon \mu_0 |X| \delta_0 /2 \le |X_i| \delta_0/2.
$$
    \item For any $j$ the set $X_S = \bigcup_{i \in S_j} X_i$ will be $(K_S, g(K_S), \delta_j/2)$-tubular because the number of elements removed from $X_S$ at steps $j' > j$ is at most
$$
2^{d+1}|X| \sum_{j' > j} \delta_{j'} \le \varepsilon \mu_0 \delta_0 2^{-(d+m+4)j} |X| \sum_{k = 1}^{\infty} 2^{-(d+m+4)k} \le \varepsilon \mu_0 \delta_0 2^{-(d+m+3)j} |X| 2^{-d-2-m}/2 = \delta_j/2. 
$$
\end{enumerate}

We clearly have $\delta_{2^m-1} \gg_{\varepsilon, \mu_0, \delta_0, m} 1$. But $m \le 1/\mu_0$ and so $\delta \gg_{K, \varepsilon, d} 1$. The lemma is proved.

\end{proofof}

\smallskip

\subsection{From tubes to subset sums}

\smallskip

The main auxiliary result in this section is the following Proposition inspired by the ideas from \cite[Section 7.2]{Zak20}. 

\begin{proposition}\label{exp}
Let $d \ge 1$, $K \ge 1$, $\delta > 0$ and $\mu > 0$ and let $K_2 > K$ be sufficiently large with respect to parameters $K, d, \delta, \mu$. Let $p > p_0(d, K, \delta, \mu)$ be a sufficiently large prime. 

Fix some $l \in [0, d]$ and let $Y \subset [-K, K]^l$ be a non-empty set. For $y \in Y$ let $X_y \subset \{y\} \times \F_p^{d-l}$ be an arbitrary (multi-)set of size at least $\mu p$. Denote $X = \bigcup_{y \in Y} X_y$.  
Suppose that $X$ is $(K_2, \delta)$-thick along any linear function which is not constant on $\{0\} \times \F_p^{d-l}$. Then there is some $u_0 \in \F_p^l$ and $k \in \N$ such that for every $u \in \F_p^{d-l}$ there are subsets $S_y \subset X_y$ such that:
\begin{align*}
    \sum_{y \in Y} \sum_{x \in S_y} x = (u_0, u),\label{s}
\end{align*}
and $\sum_{y \in Y} |S_y| = k$.
\end{proposition}

We present the proof of Proposition \ref{exp} below. The argument is based on the following lemma which was essentially proved by Alon and Dubiner \cite[Corollary 2.3 and Proposition 2.4]{AD95}. See also \cite[Lemmas 3.1 and 3.2]{Zak20}.

\begin{lemma}\label{ad}
Let $A \subset \F_p^d$ be a multiset which is $(K, \delta)$-thick along any linear function $\xi$ without constant term, for some $K \ge 0$ and $\delta > 0$. Then for any set $Y \subset \F_p^d$ of size at most $p^d/2$ there is $a \in A$ such that
\begin{equation*}
    |(Y + a) \setminus Y| \ge \max \left\{\frac{|Y|^{\frac{d-1}{d}}}{2d}, \,\frac{K\delta}{c_0 p}|Y|\right\},
\end{equation*}
where $c_0 \le 10^{10}$ is an absolute constant.
\end{lemma}

\begin{proof}
For the first inequality, we proceed like in the proof of \cite[Lemma 3.3]{Zak20} and make use of the so-called Loomis--Whitney inequality \cite{LW49}.

\begin{lemma} \label{LW} 
Let $X \subset \mathbb R^d$ be a finite set. Let $X_i$ be the projection of $X$ on the $i$-th coordinate hyperplane $\{(x_1, \ldots, x_d)~|~x_i=0\}$. Then one has an inequality $|X|^{d-1} \le \prod_{i = 1}^d |X_i|$.
\end{lemma}

Let $Y \subset \F_p^d$ and $|Y| = x^d \le p^d/2$. Since $A$ is $(K, \delta)$-thick with $\delta > 0$, there are vectors $e_1, \ldots, e_d \in A$ which form a basis of the space $\F_p^d$. 
Consider the standard embedding of $\F_p^d$ in $\mathbb Z^d$ with respect to the basis $e_1, \ldots, e_d$. Lemma \ref{LW} applied to the image of $Y$ yields that there is $i \in \{1, \ldots, d\}$ such that $|Y_i| \ge x^{d-1}$. This means that at least $x^{d-1}$ lines of the form $l_v = \{v + te_i\} \subset \F_p^d$ intersect $Y$. For any line $l_v$ intersecting $Y$ we have either $|(Y \cup (Y+e_i)) \cap l_v| > |Y \cap l_v|$ or $l_v \subset Y$. But the number of the latter lines is at most $|Y|/p$ since these lines are disjoint and contain $p$ points each. Thus, since $x \le 2^{-1/d}p$,
$$
|(Y+e_i) \setminus Y| \ge x^{d-1} - x^d/p \ge \frac{x^{d-1}}{2d}.
$$

The second inequality follows from \cite[Lemma 3.1]{Zak20}.
\end{proof}

\begin{proofof}{Proposition \ref{exp}}
Let us first consider the case $l = 0$. So, $X$ is $(K_2, \delta)$-thick along any non-constant linear function and $|X| \ge \mu p$. Then, the multiset $A = X - X$ is $(K_2, \delta)$-thick along any linear function without constant term. Indeed, suppose that for some linear function $\xi$ more than $(1-\delta)|X|^2$ differences $(x_1 - x_2)$ belong to $H(\xi, K_2)$. Then, by the pigeonhole principle, there is $x_2 \in X$ such that more than $(1-\delta)|X|$ vectors $x_1 \in X$ belong to $x_2 + H(\xi, K_2)$. But this contradicts the assumption that $X$ is $(K_2, \delta)$-thick.

Now Lemma \ref{ad} can be applied to the multiset $A$. By choosing $K_2$ sufficiently large and applying Lemma \ref{ad} iteratively one can construct a sequence of pairwise disjoint sets $\{a_1, b_1\}, \ldots, \{a_t, b_t\} \subset X$ such that
\begin{equation*} \label{CD}
    \{a_1, b_1\} + \ldots + \{a_t, b_t\} = \F_p^d.
\end{equation*}
To see this, note that at each step we can apply Lemma \ref{ad} with $Y_i = \{a_1', b_1'\} + \ldots + \{ a_i', b_i'\}$. By the thickness of $X$, if $i \le 0.1\mu \delta p$ then the set $X\setminus \left\{a_{1}',b_{1}',\ldots,a_{i}',b_{i}'\right\}$ is $(K, \delta/2)$-thick along any linear function without constant term.
In this case, by Lemma \ref{ad} we have some $x_1 - x_2 \in X-X$ so that 
$$|(Y_i+x_{1}-x_{2})\setminus Y_i| \geq \max\left\{\frac{|Y_i|^{\frac{d-1}{d}}}{2d}, \,\frac{K_2\delta}{2c_0 p}|Y_i|\right\}.$$
Since $(Y_{i}+x_{1}) \cup (Y_{i}+x_{2})=Y_{i} \cup (Y_{i} + x_{1}-x_{2})+x_{2}$, we can thus choose $\left\{a_{i+1}',b_{i+1}'\right\}=\left\{x_{1},x_{2}\right\}$ such that $Y_{i+1}:=Y_{i} + \left\{a_{i+1}',b_{i+1}'\right\}$ is significantly larger than $Y_i$. It can be easily checked that if we take $K_2$ sufficiently large then this process stops after at most $r:=0.01\delta \mu p$ steps, when we eventually reach a set $Y_{r}=\{a_1', b_1'\} + \ldots + \{a_r', b_r'\}$ of size $>p^{d}/2$ (unless we cover the whole space $\mathbb{F}_{p}^{d}$). On the other hand, if we let $X'$ denote the set $X \setminus \left\{a_1',b_1',\ldots,a_r',b_r'\right\}$, then $X'$ is also $(K_{2},0.9\delta)$-thick along any linear function without constant term since $|X \setminus X'| \leq 0.02\delta \mu p \leq 0.1\delta|X|$. Moreover, by a similar argument as above, the difference set $X'-X'$ is also $(K_{2},0.9\delta)$-thick along any linear function without constant term, so one can repeat the procedure to produce a sequence of disjoint pairs $\left\{a_1'',b_1''\right\},\ldots,\left\{a_{s}'',b_{s}''\right\} \subset X'$ such that $Y'_{s}:=\{a_1'', b_1''\} + \ldots + \{a_s'', b_s''\}$ has size $>p^{d}/2$. 
By applying the easy case of the Cauchy-Davenport theorem in $\mathbb{F}_{p}^{d}$ (for a reference, see for example \cite{EK98}), we thus have that $Y_{r}+Y'_{s} = \mathbb{F}_{p}^{d}$, so one can just take $t=r+s$ and consider $\{a_1, b_1\}, \ldots, \{a_t, b_t\}$ to be the concatenation of the two disjoint lists $\{a_1', b_1'\}, \ldots, \{a_r', b_r'\}$ and $\{a_1'', b_1''\}, \ldots, \{a_s'', b_s''\}$ so that \eqref{CD} is satisfied. Using this construction, the conclusion of Proposition \ref{exp} immediately follows with $k = t$.

Now we consider the general case, that is, $l \in [d]$ is arbitrary. In this case the multiset $X-X \subset \F_p^d$ is not necessarily $(K_2, \delta)$-thick and so we cannot apply Lemma \ref{ad}. Instead, we are going to construct a certain function $f: \{0\} \times \F_p^{d-l} \rightarrow \mathbb Q_+$ which will be $(K', \delta')$-thick along any linear function without constant term. Then we will apply Lemma \ref{ad} to $f$ in a similar manner as in the case $l = 0$ to conclude the proof. 

To define $f$ let us consider the set $\Lambda \subset \mathbb Z^{Y}$ consisting of all integer vectors $(\lambda_y)_{y \in Y}$ such that $\sum \lambda_y y = 0$ and $\sum \lambda_y = 0$.
For a pair of vectors $\lambda^1, \lambda^2 \in \N^{Y}$ let $\mathcal J^{\lambda^1, \lambda^2}$ be the set of all pairs $(J_1, J_2)$, where $J_1$ and $J_2$ are disjoint labeled subsets of $X$ and for every $y \in Y$ we have 
\begin{align*}
    |J_1 \cap X_y| = \lambda_y^1, \\
    |J_2 \cap X_y| = \lambda_y^2,
\end{align*}
For a pair of (multi-)sets $(J_1, J_2)$ we denote by $\sigma(J_1, J_2)$ the sum of elements of $J_1$ minus the sum of elements of $J_2$. Note that if $\lambda^1 - \lambda^2 \in \Lambda$ then for any $(J_1, J_2) \in \mathcal J^{\lambda^1, \lambda^2}$ we have $\sigma(J_1, J_2) \in \{0\} \times \F_p^{d-l}$ and $|J_1| = |J_2|$.

For $\lambda^1, \lambda^2 \in \N^Y$ such that $\lambda^1 - \lambda^2 \in \Lambda$ let $f_{\lambda^1, \lambda^2}: \{0\} \times \F_p^{d-l} \rightarrow \N$ be the characteristic function of the multiset of all sums $\sigma(J_1, J_2)$ over $(J_1, J_2) \in \mathcal J^{\lambda^1, \lambda^2}$. We define the function $f: \{0\} \times \F_p^{d-l} \rightarrow \mathbb Q_+$ as follows:
\begin{equation}\label{f}
    f = \sum_{\lambda^1, \lambda^2} \frac{f_{\lambda^1, \lambda^2}}{|\mathcal J^{\lambda^1, \lambda^2}|},
\end{equation}
where the sum is taken over all vectors $\lambda^1, \lambda^2 \in \N^Y$ with $l_1$-norm bounded by a sufficiently large constant $T = T(d, K, \delta, \mu) > 0$ and such that $\lambda^1 - \lambda^2 \in \Lambda$.

Note that we have the following expression for the size of $\mathcal J^{\lambda^1, \lambda^2}$, as a product of multinomial coefficients:
\begin{equation}\label{J}
    |\mathcal J^{\lambda^1, \lambda^2}| = \prod_{y \in Y} {|X_y| \choose \lambda^1_y, \lambda^2_y}.
\end{equation}

It requires some work to show that $f$ is indeed $(K', \delta')$-thick along any linear function on $\{0\} \times \F_p^{d-l}$, so we isolate this fact as a separate lemma. 

\begin{lemma} \label{Athick}
The function $f$ is $(K', \delta')$-thick along any linear function without constant term on $\{0\} \times \F_p^{d-l}$. Here $K'$ and $\delta'$ depend on parameters $K, K_2, d, \delta, \mu$ in such a way that $K'$ can be arbitrarily large compared to $K, d, \delta, \mu$ and $\delta'$ if one takes $K_2$ large enough.
\end{lemma}

A somewhat similar result was proved in \cite[Lemma 7.3]{Zak20} but the details are different, so we include below the complete proof for our setup.

\bigskip

{\textbf{Proof of Lemma \ref{Athick}:}} Suppose that $f$ is not $(K', \delta')$-thick along some linear function $\xi\ :\ \{0\} \times \mathbb{F}_{p}^{d-l} \to \mathbb{F}_{p}$ without constant term. Here we take $\delta' = \min\{0.01/T, 0.1\delta\}$. The number $K'$ will be determined at the end of the argument.

We may extend $\xi$ to a linear function on $\F_p^d$ which does not depend on the first $l$ coordinates. 
Denote by $S \subset \mathbb Z^Y \times \mathbb Z^Y$ the set of all pairs $\lambda^1, \lambda^2 \in \N^Y$ satisfying $\|\lambda^i\|_1 \le T$ and $\lambda^1 - \lambda^2 \in \Lambda$. Then by (\ref{f}) we have $f(\{0\} \times \mathbb F^{d-l}_p) = |S|$.
So we have $f(H(\xi, K')) \ge (1-\delta') |S|$. 
Let $S' \subset S$ be the set of all pairs $(\lambda^1, \lambda^2) \in S$ such that $f_{\lambda^1, \lambda^2}(H(\xi, K')) \ge (1 - 2\delta') |\mathcal J^{\lambda^1, \lambda^2}|$, that is, the function $f_{\lambda^1, \lambda^2}$ is not $(K', 2\delta')$-thick along $\xi$. Note that
$$
\delta' |S| \ge f(\mathbb F_p^d \setminus H(\xi, K)) \ge \sum_{(\lambda^1, \lambda^2) \in S\setminus S'} \frac{f_{\lambda^1, \lambda^2}(\mathbb F_p^d \setminus H(\xi, K))}{|\mathcal J^{\lambda^1, \lambda^2}|} \ge 2\delta' |S \setminus S'|,
$$
and so $|S'| \ge |S| /2$.

\begin{claim}\label{cl1}
For any $y \in Y$ there is $r_y \in \mathbb F_p$ such that $|X_y \cap H(\xi - r_y, 2K')| \ge (1-10\delta') |X_y|$. 
\end{claim}

\begin{proof}
Recall that $\Lambda$ is a lattice defined by equations with coefficients bounded by $K$. Thus, if $T$ is large enough compared to $K$ then any subspace $V \subset \mathbb R^Y \times \mathbb R^Y$ either contains $S$ or intersects it in at most $|S|/10$ elements. In particular, for any subspace $V$ we have $S' \not\subset V$ or $S \subset V$.

For $y \in Y$ consider the subspace $V = \{(\lambda^1, \lambda^2) \in \mathbb R^Y \times \mathbb R^Y ~|~ \lambda^1_y = \lambda^2_y = 0\}$. Then we clearly have $S \not \subset V$ because the vector $(e_y, e_y)$, where $e_y$ has only one non-zero coordinate in $y$-th place, belongs to $S$. So there is a pair $(\lambda^1, \lambda^2) \in S'$ such that $\lambda^i_y \neq 0$ for some $i \in \{1, 2\}$. Without loss of generality we may assume that $i =1$.

For labeled elements $x_1, x_2$ of $X_y$ denote by $\mathcal J^{\lambda^1, \lambda^2}_{x_1, x_2}$ the subfamily in $\mathcal J^{\lambda^1, \lambda^2}$ consisting of all pairs $(J_1, J_2)$ such that $x_1 \in J_1$, $x_2 \not \in J_1$ and $x_1, x_2 \not \in J_2$. Define a graph $G$ on the set of labeled elements of $X_y$ where a pair $x_1, x_2$ of labeled elements of $X_y$ forms an edge if $\xi(x_1 - x_2) \not \in [-2K', 2K']$. Let $x_1, x_2$ be an edge and take $(J_1, J_2) \in \mathcal J_{x_1, x_2}^{\lambda^1, \lambda^2}$. Denote $J_1' = J_1 \setminus \{x_1\} \cup \{x_2\}$, then we have $(J_1', J_2) \in \mathcal J_{x_2, x_1}^{\lambda^1, \lambda^2}$ and 
$$
\sigma(J_1, J_2) - \sigma(J_1', J_2) = x_1 - x_2.
$$
Thus, one of the elements $\xi(\sigma(J_1, J_2))$ or $\xi(\sigma(J_1', J_2))$ does not belong to the interval $[-K', K']$. So at least half of the sums $\sigma(J_1, J_2)$, $\sigma(J'_1, J_2)$ over $(J_1, J_2) \in \mathcal J_{x_1, x_2}^{\lambda^1, \lambda^2}$ does not belong to the strip $H(\xi, K')$. We are going to combine this observation with the fact that the multiset of sums $\sigma(J_1, J_2)$ over $(J_1, J_2) \in \mathcal J^{\lambda^1, \lambda^2}$ is not $(K', \delta')$-thick along $\xi$ to conclude that the graph $G$ has a very large independence number. 

Indeed, suppose that the independence number of $G$ is at most $|X_y| - r$ for some $r \ge 1$. Then one can find $r$ disjoint edges $\{x_1, y_1\}, \ldots, \{x_r, y_r\}$ in the graph $G$. By an elementary double counting argument, for any $i \neq j$, $i, j \in [r]$, we have 
$$
|\mathcal J_{x_i, y_i}^{\lambda^1, \lambda^2} \cap \mathcal J_{x_j, y_j}^{\lambda^1, \lambda^2}| = \frac{(\lambda^1_y - 1)(|X_y| - \lambda^1_y - \lambda^2_y - 1)}{(|X_y| - 2)(|X_y| - 3)} |\mathcal J_{x_i, y_i}^{\lambda^1, \lambda^2}| \le \frac{\lambda^1_y}{|X_y|} |\mathcal J_{x_i, y_i}^{\lambda^1, \lambda^2}|.
$$
Applying Bonferroni's inequality to the sets $\mathcal J^{\lambda^1, \lambda^2}_{x_i, y_i} \cup \mathcal J^{\lambda^1, \lambda^2}_{y_i, x_i}$, $i = 1, \ldots, r'$, for some $r' \le r$,
we can estimate the number $N$ of pairs $(J_1, J_2) \in \mathcal J^{\lambda^1, \lambda^2}$ such that $\xi(\sigma(J_1, J_2)) \not \in [-K', K']$:
\begin{align*}
    N \ge |\mathcal J^{\lambda^1, \lambda^2}_{x_1, y_1}|\left(r' - 4{r' \choose 2} \frac{\lambda^1_y}{|X_y|}\right).
\end{align*}
A double counting argument shows that 
$$
|\mathcal J^{\lambda^1_y, \lambda^2_y}_{x_1, y_1}| = \frac{\lambda^1_y (|X_y| - \lambda^1_y - \lambda^2_y)}{|X_y|(|X_y| - 1)} |\mathcal J^{\lambda^1, \lambda^2}| \ge \frac{\lambda^1_y}{2|X_y|}|\mathcal J^{\lambda^1, \lambda^2}|,
$$
since $|X_y| \ge \mu p$, $\lambda^i_y \le \|\lambda^i\|_1 \le T$ and $p$ is large enough compared to $\mu$ and $T$. We conclude that for any $r' \le r$:
$$
N \ge \frac{\lambda^1_y}{2|X_y|}\left (r' - 2r'^2 \frac{\lambda^1}{|X_y|}\right)|\mathcal J^{\lambda^1, \lambda^2}|.
$$
On the other hand, since $(\lambda^1, \lambda^2) \in S'$ we have $N \le 2\delta' |\mathcal J^{\lambda^1, \lambda^2}|$. It is easy to see that the last two inequalities are incompatible if we put $r' = 10\delta' \frac{|X_y|}{\lambda^1_y}$ and note that $\delta' \le 0.01$.

Thus, the independence number of the graph $G$ is at least $(1 - 10\delta')|X_y|$. Let $Z'_y$ be an independent set of labeled elements of $X_y$ of maximal size. Pick an arbitrary $z \in Z'_y$ and  put $r_y = \xi(z) \in \F_p$. Then by the definition of $G$ we have 
$$
Z'_y \subset H(\xi - r_y, 2K'),
$$
and $|Z'_y| \ge (1 - 10\delta') |X_y|$ and the claim follows.
\end{proof}

Denote $Z_y = X_y \cap H(\xi-r_y, 2K')$ so that we have $|Z_y| \ge (1 - 10\delta') |X_y|$. The next step is to show that the value $r_y$ can be approximated by a linear function of $y \in Y$. This will then contradict the thickness condition on the multiset $X$.


For a given $(\lambda^1, \lambda^2) \in S$ denote by $\mathcal J^{\lambda^1, \lambda^2}_0$ a subfamily in $\mathcal J^{\lambda^1, \lambda^2}$ consisting of those pairs $(J_1, J_2)$ such that for every $(x, n) \in J_1 \cup J_2$ we have $x \in \bigcup_{y \in Y} Z_y$. From (\ref{J}) and a similar expression for the size of $\mathcal J^{\lambda^1, \lambda^2}_0$ we see that
$$
\frac{|\mathcal J^{\lambda^1, \lambda^2}_0|}{|\mathcal J^{\lambda^1, \lambda^2}|} \ge 1 - 20\delta' (\|\lambda^1\|_1 + \|\lambda^2\|_1) \ge 1 - 40 \delta' T.
$$
Here we used the fact that $|X_y| \ge \mu p$ and $p$ is large enough which implies that ${|X_y| \choose a, b} = (1 + o_p(1)) \frac{|X_y|^{a+b}}{a! b!}$ for $a, b \le T$. Since $\delta' \le 0.01/ T$, if $(\lambda^1, \lambda^2) \in S'$ then there is a pair $(J_1, J_2) \in \mathcal J^{\lambda^1, \lambda^2}_0$ such that $\sigma(J_1, J_2) \in H(\xi, K')$. 

But $|r_y - x| \le 2K'$ for every $x \in Z_y$. And so, slightly abusing notation, we can estimate:
\begin{align*}
    \left| \sum_{y \in Y} r_y (\lambda^1_y - \lambda^2_y)\right| \le 4K'T + \left | \sum_{(x, n) \in J_1} \xi(x) - \sum_{(x, n) \in J_2} \xi(x) \right| \le 4K' T + K' \le 5K'T.
\end{align*}
In other words, for any $(\lambda^1, \lambda^2) \in S'$ we have $\langle r, \lambda^1 - \lambda^2\rangle \in [-5K'T, 5K'T]$ where $\langle \cdot, \cdot \rangle$ denotes the standard pairing on $\F_p^Y$. Using the fact that $|S'| \ge |S|$ and some elementary linear algebra one can find a constant $C\ll_{K'T, |Y|} 1$ such that for any vector $\lambda \in \Lambda$ with $l_1$-norm at most $T$ we have $\langle r, \lambda \rangle \in [-C, C]$.

Let $Y' \subset Y$ be a minimal set whose affine hull contains $Y$. Then for every $y \in Y \setminus Y'$ there is a unique up to scaling vector $\lambda(y) \in \Lambda$ such that $\lambda(y)_y \neq 0$ and $\lambda(y)_{\tilde y} = 0$ for any $\tilde y \not \in Y' \cup \{y\}$. We fix $\lambda(y)$ in such a way that $\lambda(y)_y$ has minimal positive value. Note that if $p$ is large enough then the number $\lambda(y)_y$ is not divisible by $p$. Now we define a new vector $\tilde r = (\tilde r_y) \in \F_p^Y$ as follows. For $y \in Y'$ we let $\tilde r_{y} = r_y$ and for $y \in Y \setminus Y'$ we define
$$
\tilde r_y = - \frac{1}{\lambda(y)_y}\sum_{y' \in Y'} \lambda(y)_{y'} r_{y'}.
$$
It is now clear that $\lambda(y)_y \tilde r_y - \lambda(y)_y r_y \in [-C, C]$ (considered as an element of $\F_p$). Let $\eta: \F_p^l \rightarrow \F_p$ be a linear function such that $\eta(y) = r_y$ for $y \in Y'$. Then from linearity we have $\eta(y) = \tilde r_y$ for every $y \in Y$. Denote $L = \prod_{y \in Y\setminus Y'} \lambda(y)_y$. We may view $\eta$ as a function on $\F_p^d$ which does not depend on the last $d-l$ coordinates. Then for any $x \in Z_y$ we have 
$$
|L\xi(x) - L\eta(x)| \le 2K'L + |L r_y - L \tilde r_y| \le 2K'L + LC =: B,
$$
i.e. if we let $\xi' = L\xi - L\eta$ then for any $x \in Z_y$ we have $\xi'(x) \in [-B, B]$. Therefore, the multiset $X$ is not $(B, 10\delta')$-thick along $\xi'$. It is easy to see that $B$ is bounded in terms of parameters $K', K, |Y|$ and $d$. Define $K'$ to be the largest number such that $B \le K_2$. Then it is clear that $K'$ can be made arbitrarily large compared to $K, \delta$ and $\delta'$ if one takes $K_2$ sufficiently large. Indeed, recall that $\delta' = \min\{0.01/T, 0.1\delta\}$ where $T$ depends on $K$ and $d$ only. With these parameters, the multiset $X$ is not $(K_2, \delta)$-thick along the linear function $\xi'$. But $\xi'$ coincides with $L\xi$ on $\{0\} \times \F_p^{d-l}$ and therefore is not constant on this subspace. This contradicts to the initial assumption on $X$.

$\hfill \square$

\bigskip

Returning to the proof of Proposition \ref{exp}, we can now similarly apply Lemma \ref{ad} repeatedly to the multiset $A$ like in the case $l=0$, and construct recursively a sequence of pairs $\left(J_1^1, J_2^1\right), \ldots, \left(J_1^t, J_2^t\right)$ where all sets $J_i^j$ are pairwise disjoint as labeled subsets of $X$ and such that
\begin{equation*} \label{CD2}
    \{\sigma(J_1^1), \sigma(J_2^1)\} + \ldots + \{\sigma(J_1^t), \sigma(J_2^t)\} = \{u_0\} \times \F_p^{d-l},
\end{equation*}
for some fixed vector $u_0$ (by first producing two pairwise disjoint lists $\left(J_{1,1}^1, J_{1,2}^1\right), \ldots, \left(J_{1,1}^{r}, J_{1,2}^{r}\right)$ and $\left(J_{2,1}^1, J_{2,2}^1\right), \ldots, \left(J_{2,1}^{s}, J_{2,2}^{s}\right)$ such that
$\left|\{\sigma(J_{1,1}^1), \sigma(J_{1,2}^1)\} + \ldots + \{\sigma(J_{1,1}^r), \sigma(J_{1,2}^r)\}\right| > p^{d-l}/2$ and $\left|\{\sigma(J_{2,1}^1), \sigma(J_{2,2}^1)\} + \ldots + \{\sigma(J_{2,1}^s), \sigma(J_{2,2}^s)\}\right| > p^{d-l}/2$, and then using again the Cauchy-Davenport theorem in $\{u_0\} \times \F_p^{d-l}$ to see that the concatenation of the two lists (of length $t=r+s$) satisfies \eqref{CD2}. Using this construction, Proposition \ref{exp} immediately follows with $k$ equal to $|J_1^1| + \ldots + |J_1^t|$ (note that by definition $|J_1^i| = |J_2^i|$). 
\end{proofof}

\smallskip

\subsection{High multiplicity case}

\smallskip

The final lemma is a result about zero sums in sequences which can be regarded as a generalization of Olson's main result from \cite{Ols69}.

\begin{lemma}\label{nul}
Let $Y \subset \F_p^d$ be an arbitrary set and let $w: Y \rightarrow \N$ be a function such that $\sum_{y \in Y}w(y) \ge d(p-1) + 2r|Y| + 1$ for some $r \ge 0$. Then, there exist coefficients $a_y \in \mathbb{N}$, one for each $y \in Y$, such that $a_y \in \{0\} \cup [r, w(y)-r]$ and $\sum_{y \in Y} a_y y = 0$, while not all $a_{y}$'s are simultaneously zero.
\end{lemma}

When $r=0$, notice that this indeed immediately implies that if $n > d(p-1)$, then among any $n$ elements $v_{1},\ldots,v_{n}$ of $\mathbb{F}_{p}^{d}$ there exists a nonempty subsequence with a zero subsum. To prove Lemma \ref{nul}, we will make use of Alon's Combinatorial Nullstellensatz \cite[Theorem 1.2]{NA99}, which we recall for the reader's convenience. 

\begin{lemma}\label{CN}
Let $\mathbb{F}$ be an arbitrary field, and let $f = f(x_{1},\ldots,x_{n})$ be a polynomial in $\mathbb{F}[x_{1},\ldots,x_{n}]$. Suppose the degree $\deg(f)$ of $f$ is $\sum_{i=1}^{n}{t_{i}}$, where each $t_{i}$ is a nonnegative integer, and suppose that the coefficient of $\prod_{i=1}^{n}{x_{i\
}^{t_{i}}}$ in $f$ is nonzero. Then, if $S_{1},\ldots,S_{n}$ are the subsets of $F$ with $|S_{i}| > t_{i}$, there exist $s_{1} \in S_{1}$,$\ldots$, $s_{n} \in S_{n}$ so that
$$f(s_{1},\ldots,s_{n}) \neq 0.$$
\end{lemma}

\begin{proofof}{Lemma \ref{nul}}
For $y \in Y$ denote $A_y = \{0\} \cup [r, w(y)-r]$ and consider the following $|Y|$-variate polynomial in $(\alpha_y)_{y \in Y}$:
$$
P(\alpha_y~|~ y\in Y) = \prod_{i=1}^d \left( 1 - \left( \sum_{y \in Y}  \alpha_y y_i\right)^{p-1} \right),
$$
where $y_i$ denotes the $i$-th coordinate of $y$ as an element in $\mathbb{F}_{p}^{d}$. Note that $P(\alpha_y~|~ y\in Y)$ is non-zero if and only if $\sum_{y\in Y} \alpha_y y =0$, so the zero vector $\overrightarrow{0}$ in $\mathbb{F}_{p}^{|Y|}$ is certainly not a zero of the polynomial $P$. On the other hand, observe that
$$
\sum_{y \in Y} (|A_y| - 1) \ge \sum_{y \in Y} w(y) - 2r|Y| > d(p-1),
$$
so by Lemma \ref{CN} applied in a slightly smaller cartesian product which is strictly contained in $\prod_{y \in Y} A_y$ and which does not contain $\overrightarrow{0}$, it follows that $P$ must take some other non-zero value at a vector in $\prod_{y \in Y} A_y$ that does not have all coordinates equal to $0$. This completes the proof of Lemma \ref{nul}.

\end{proofof}

\bigskip

\section{Proof of Theorem \ref{main}}

\smallskip

Let $X \subset \F_p^d$ be an arbitrary set of size $(d-1+\varepsilon)p$ where $p$ is a sufficiently large prime number. Let $g: \N \rightarrow \N$ be a sufficiently fast growing function.

Apply Lemma \ref{sdec} to $X$ with $\varepsilon' = \varepsilon / 2d$ and $g = g$, $K = K_0$. After removing $X_0$ from $X$ we will obtain a set $X$ of size at least $(d-1+\varepsilon/2)p$ and a decomposition $X = X_1 \cup \ldots \cup X_m$ with several important properties. In particular, if $U_i$ denotes the affine hull of the set $X_i$, recall that $X_i$ is $(g^{d+1}(K), \delta)$-thick in $U_i$ for some $\delta \gg_{K, \varepsilon, d} 1$, for each $i \in \left\{1,\ldots,m\right\}$. Moreover, for every $i$ we also have that $|X_i| \ge \mu |X|$, where $\mu \gg_{K, \varepsilon, d} 1$. Note that we may assume that $\mu$ is small enough, namely, $\mu m < \varepsilon/100$. Furthermore, note that since $X$ is a set, all spaces $U_i$ are non-zero dimensional. 

Let $H \subset \F_p^d$ be a generic hyperplane passing through the origin which intersects all affine spaces $U_i$. Such $H$ exists since $m \ll_{K, d, \varepsilon} 1$ and $p$ is large enough (indeed, a random hyperplane $H$ intersects an affine subspace of dimension $\ge 1$ with probability $\gtrsim 1- \frac{1}{p}$).


Fix an arbitrary vector $x_i \in H \cap U_i$ for each $i \in [m]$. Assign the weight $w_i = |X_i|$ and apply Lemma \ref{nul} to the set $Y = \{x_1, \ldots, x_m\} \subset H$ with the weight $w$ and $r = \mu |X|/3$. We have $\mu < \varepsilon / 10m$ and so
$$
\sum_{i = 1}^{m} w_i = |X| \ge (d-1 + \varepsilon/2)p > (d-1) p + 2\mu |X|m/3. 
$$
Thus, there are non-negative, not all zero, coefficients $a_i \in \{0\} \cup [\mu |X|/3, w_i - \mu |X|/3]$ such that $\sum_{i = 1}^m a_i x_i = 0$. Let $S \subset [m]$ be the set of $i \in [m]$ for which $a_i > 0$. By Lemma \ref{sdec}, the set $X_S = \bigcup_{i \in S} X_i$ is $(K_S, g(K_S), \mu)$-tubular for $K_S = g^l(K)$ and $l \in [0, d]$. So after a linear change of coordinates, there is a vector $v \in \F_p^l \times \{0\}$ such that 
\begin{equation*}
    X_S \subset v + [-K_S, K_S]^l \times \F_p^{d-l}.
\end{equation*}
Denote by $\pi$ the projection onto first $l$ coordinates. The condition that $X_i$ is $(g^{d+1}(K), \delta)$-thick in $U_i$ and the fact that $K_S < g^{d+1}(K)$ implies that $\pi(U_i)$ is a single point for any $i \in S$. 
Indeed, for $j = 1, \ldots, l$ consider the linear function $\xi_j(x) = x_j - v_j$. Then for any $i \in S$ the set $X_i$ is not $(K_S, 0)$-thick along $\xi_j$. This implies that $\xi_j$ is constant on $U_i$ and so $\pi(U_i) \subset \F_p^l$ is a single point.
Denote this point by $y_i \in v + [-K_S, K_S]^l$ and observe that 
\begin{equation}\label{ys}
    \sum_{i = 1}^m a_i y_i = 0,
\end{equation}
since $\pi$ is a linear operator. 

Denote by $Y \subset [-K_S, K_S]^l$ the set obtained from the projection $\pi(X)$ and by shifting by $v$.
Note that for any $y \in Y$ the set $X_y = X_S \cap \left(\{y\} \times \F_p^{d-l}\right)$ has size at least $\mu |X| \ge \mu |X_S|$ by Lemma \ref{sdec}. Moreover, the set $X_S$ is $(g(K_S), \delta)$-thick along any linear function which is non-constant on $\{0\} \times \F_p^{d-l}$. For $y \in Y$ denote by $a_y$ the sum of all numbers $a_i$ over all $i$ is such that $y = y_i - v$. 

\begin{proposition}\label{sub}
If $p$ is large enough then there are sets $Z_y \subset X_y$, such that for any $y \in Y$ we have $|Z_y| \in [\mu |X_y| / 20, \mu|X_y|/10]$ and the set $Z = \bigcup_{y\in Y} Z_y$ is $(g(K_S), \delta/4)$-thick along any linear function which is not constant on $\{0\} \times \F_p^{d-l}$.
\end{proposition}

We prove Proposition 2 by using a probabilistic argument, where we make use of the following Chernoff bound (see, for example, \cite[Corollary 21.7]{JLR00}). 



\begin{lemma} \label{chernoff}
Let $X$ be a random variable with the binomial distribution $\operatorname{Bin}(N,p)$ and $\eta \in (0,1)$. Then
\begin{eqnarray*}
\operatorname{Pr}\left(X \leq (1-\eta)\mathbb{E}[X]\right) &\leq& \exp\left(-\frac{\eta^2}{2} \mathbb{E}[X]\right)\\
\operatorname{Pr}\left(X \geq (1+\eta)\mathbb{E}[X]\right) &\leq& \exp\left(-\frac{\eta^2}{3} \mathbb{E}[X]\right)
\end{eqnarray*}
\end{lemma}

\begin{proofof}{Proposition \ref{sub}}
Choose sets $Z_y \subset X_y$ at random according to the binomial distribution $\operatorname{Bin}(|X_y|,\mu/15)$. It follows from Lemma \ref{chernoff} and the fact that $|X_y| \gg p$ that with high probability we have $|Z_y| \in [\mu |X_y| / 20, \mu|X_y|/10]$ for all $y \in Y$.
We will show that for any fixed linear function the event that $Z = \bigcup_{y \in Y}Z_y$ is not $(g(K_S), \delta/4)$-thick has probability exponentially small in $p$. 
Since there are only $O(p^d)$ linear functions on $\F_p^d$ this will be enough to prove Proposition \ref{sub}. Fix a linear function $\xi$ which is not constant on $\{0\} \times \F_p^{d-l}$ and denote $X'_y = X_y \setminus H(\xi, g(K_S))$. Since the set $X_S$ is $(g(K_S), \delta)$-thick along $\xi$, the set $X' = \bigcup_{y \in Y} X'_y$ has size at least $\delta|X_S|$. 
Note that the expected size of the intersection $Z_y \cap X'_y$ is asymptotically equal to $\frac{|Z_y| |X'_y|}{|X_y|}$ and so, provided that $|X'_y| \gg p$, by Lemma \ref{chernoff} the probability of the event that $|Z_y \cap X'_y| < \frac{|Z_y| |X'_y|}{1.5|X_y|}$ is at most $e^{-cp}$ for some $c \gg 1$. Therefore, the probability that 
$$
\sum_{y \in Y} |Z_y \cap X'_y| < \sum_{y \in Y} \frac{|Z_y| |X'_y|}{1.5|X_y|}
$$
is at most $|Y| e^{-cp}$. But since the set $Z_y$ has size in the interval $[\mu |X_y| / 20, \mu|X_y|/10]$ the right hand side is at least
$$
(\mu /20) \sum_{y \in Y} |X'_y|/1.5 \ge (\mu/20) \delta |X| /1.5 > \delta |Z| /4,
$$
which means that $Z$ is $(g(K_S), \delta/4)$-thick along $\xi$ with probability at least $1 - |Y| e^{-cp}$. This completes the proof.
\end{proofof}

Fix sets $Z_y$ as in Proposition \ref{sub}.
Now let the function $g$ grow so fast that we have $g(K) > K_2$ where $K_2 = K_2(K, d, \delta/4, \mu^2/10)$ is the function from Proposition \ref{exp}. We can then apply Proposition \ref{exp} to sets $Z_y-v$ with $u = 0$ to get some $k_y \in \N$ such that
$$
\sum_{y \in Y} k_y (y-v) = u_0,~~\sum_{y \in Y} k_y = k,
$$
where $u_0$ and $k$ are from the statement of Proposition \ref{exp}. Note that $k_y \le |Z_y| \le \mu|X_y|/10$.

For each $y \in Y$ fix a subset $A_y \subset X_y \setminus Z_y$ of size $a_y - k_y$ (which is possible thanks to the estimates on $a_i$). Let $u = (u_1, u_2) \in \F_p^d$ denote the following vector:
$$
u = \sum_{y \in Y} \sum_{x \in A_y} x.
$$
From (\ref{ys}) and from the definition of the $k_y$'s we see that, in fact, $u_1 = -u_0-kv$.

By the conclusion of Proposition \ref{exp}, applied to the vector $-u_2 \in \F_p^{d-l}$, we obtain some sets $S_y \subset Z_y$ such that $\sum_{y \in Y}|S_y| = k$ and
\begin{equation*}
    \sum_{y \in Y} \sum_{x \in S_y} (x - v) = (u_0, -u_2),
\end{equation*}
After rearranging, this rewrites as
$$
\sum_{y \in Y} \sum_{x \in S_y} x = (u_0 + k v, -u_2) = - \sum_{y \in Y}\sum_{x \in A_y} x.
$$
So we see that the set $B = \bigcup_{y \in Y} A_y \cup S_y$ has zero sum. Theorem \ref{main} is proved.


\end{document}